\documentclass[12pt]{amsart}
\usepackage{amssymb}
\usepackage{textcomp}
\usepackage[all]{xy}

\title{Basic finite \'{e}tale equivalence relations}
\author{Ying Zong}

\address{Department of Mathematics\\University of Toronto}
\email{zongying@math.utoronto.ca}
\date{} 

\begin{document}

\maketitle


1. \emph{Introduction.}
\smallskip

Let $S$ be an algebraic space, $A$ an $S$-abelian algebraic space, $X$ an $A$-torsor on $S$ for the \'{e}tale topology and $L$ a finite \'{e}tale $S$-equivalence relation on $X$. Thus, $X$ is representable by an $S$-algebraic space and its quotient, in the sense of \'{e}tale topology, by the finite \'{e}tale $S$-equivalence relation $L$, is representable by a proper and smooth $S$-algebraic space $Y$ with geometrically irreducible fibers. In the following we define the property for $L$ to be a \emph{basic} finite \'{e}tale equivalence relation and show that this property is remarkably rigid and our first goal is to prove :
\smallskip

{\bf Theorem 1.1.} --- \emph{Let $f: Y\to S$ be a proper smooth morphism of algebraic spaces with $S$ connected. Assume that for one geometric point $t$ of $S$ there is a finite \'{e}tale $t$-morphism from a $t$-abelian variety $A'_t$ onto $Y_t=f^{-1}(t)$ with degree $[A'_t: Y_t]$ prime to the residue characteristics of $S$.}
\smallskip

\emph{Then }:
\smallskip

\emph{a) There exist an $S$-abelian algebraic space $A$, an $A$-torsor $X$ on $S$ for the \'{e}tale topology and a basic finite \'{e}tale $S$-equivalence relation $L$ on $X$ such that $Y$ is $S$-isomorphic to the quotient $X/L$.}
\smallskip

\emph{b) Let $(A, X, L)$ be as in a). Then $X$ is a $f^*G$-torsor on $Y$ for the \'{e}tale topology where the $S$-group $G=f_*\underline{\mathrm{Aut}}_Y(X)$ is a finite \'{e}tale $S$-algebraic space. The functor 
\[(A', X', L')\mapsto P'=f_*\underline{\mathrm{Isom}}_Y(X, X')\] induces a bijection from the set of all triples as in a) up to isomorphisms onto $H^1(S, G)$.}

\smallskip

The proof is in \S 9. 
\smallskip

When $f$ degenerates to $\overline{f}: \overline{Y}\to \overline{S}$ where $\overline{f}$ is separated of finite type and universally open with $\overline{S}$ locally noetherian, then a geometric fiber $\overline{f}^{-1}(\overline{s})$, if it is non-empty and does not have uniruled irreducible components, is irreducible; see \S 10, which, as a continuation of \cite{almost} and our main objective, contains further support for the \emph{Principle} : \emph{Non-uniruled abelian degenerations are almost non-degenerate}.
\smallskip

2. \emph{Definition.}
\smallskip

{\bf Definition 2.1.} --- \emph{Let $S$ be an algebraic space, $A$ an $S$-abelian algebraic space, $X$ an $A$-torsor on $S$ for the \'{e}tale topology and $L$ a finite \'{e}tale $S$-equivalence relation on $X$. We say that $L$ is basic at a geometric point $s$ of $S$ provided the following condition holds }:
\smallskip

\emph{If $E$ is a finite \'{e}tale $s$-subgroup of $A_s$ such that the morphism
\[E\times_s X_s \to X_s\times_sX_s,\ (e, x)\mapsto (e+x, x)\] factors through the graph of $L_s$, then $E=0$.}
\smallskip

\emph{We say that $L$ is basic if it is basic at every geometric point of $S$.}
\smallskip

Thus, put in words, $L$ is basic at $s$ if and only if no sub-equivalence relation of $L(s)$ other than the identity is generated by a translation. Let $Y$ denote the quotient $X/L$. The graph of $L$ is then identical to $X\times_YX$. With the above notations, saying that the morphism
\[E\times_sX_s\to X_s\times_sX_s,\ (e, x)\mapsto (e+x, x)\] factors through the graph of $L_s$ amounts to saying that the translations of $X_s$ defined by the elements of $E(s)$ are $Y_s$-automorphisms.
\smallskip

2.2. Let $s'\to s$ be a morphism of geometric points of $S$. The base change
\[L_s\mapsto L_s\times_ss'\] establishes a bijection between the collection of finite \'{e}tale $s$-equivalence relations on $X_s$ which are basic at $s$ and the collection of finite \'{e}tale $s'$-equivalence relations on $X_{s'}$ which are basic at $s'$. Indeed, for every proper $s$-algebraic space $Z$, the functor of base change by $s'\to s$ is an equivalence from the category of finite \'{e}tale $Z$-algebraic spaces onto the category of finite \'{e}tale $Z_{s'}$-algebraic spaces (SGA 4 XII 5.4). This same base change provides also an equivalence between the category of finite \'{e}tale $s$-subgroups of $A_s$ and the category of finite \'{e}tale $s'$-subgroups of $A_{s'}$.
In particular, a given finite \'{e}tale $S$-equivalence relation $L$ on $X$ is basic at $s$ if and only if it is basic at $s'$. The collection of basic finite \'{e}tale equivalence relations on $X_s$ and on $X_{s'}$ are in canonical bijective correspondence.
One deduces also that if $L$ is basic, so is $L_{S'}$ for every base change $S'\to S$. When $S'\to S$ is surjective, the condition that $L$ be basic is then equivalent to that $L_{S'}$ be basic.

\smallskip

3. \emph{Automorphisms I.}
\smallskip

We study now the $Y$-automorphisms of $X$. 
\smallskip

{\bf Lemma 3.1.} --- \emph{Keep the notations of $(2.1)$. With each $S$-morphism $q: X\to X$ there is associated a unique $S$-group homomorphism $p: A\to A$ such that $q$ is $p$-equivariant.}

\begin{proof} One may assume the torsor $X$ to be trivial, as the existence of a unique homomorphism $p$ is an \'{e}tale local question on $S$.
\smallskip

Each $S$-morphism $q: A\to A$ is the unique composite of a translation ($x\mapsto x+q(0)$) and an $S$-group homomorphism $p: A\to A$ (``Geometric Invariant Theory'' 6.4). It is clear that $q$ is $p$-equivariant : for every two local $S$-sections $a, x$ of $A$, one has
\[q(a+x)=p(a)+q(x).\]

\end{proof}
\smallskip

{\bf Lemma 3.2.} --- \emph{Keep the notations of $(2.1)$. Assume that $L$ is basic. Then every $Y$-endomorphism of $X$ is a $Y$-automorphism. If $X$ admits furthermore an $S$-section $x$, the only $Y$-automorphism of $X$ preserving $x$ is the identity morphism.}

\begin{proof} Every $Y$-morphism $X\to X$ is finite \'{e}tale, as $X$ is finite \'{e}tale over $Y$. Let $q$ be such a $Y$-endomorphism of $X$. By (3.1), $q$ is equivariant relative to a unique $S$-group homomorphism $p: A\to A$. And, $p$ is finite \'{e}tale, since it is \'{e}tale locally on $S$ isomorphic to $q$. In particular, $E=\mathrm{Ker}(p)$ is a finite \'{e}tale $S$-subgroup of $A$. Observe that the morphism
\[E\times_SX\to X\times_SX,\ (e, x)\mapsto (e+x, x)\] factors through the graph of $L$. That is to say, for every local $S$-section $e$ (resp. $x$) of $E$ (resp. $X$), $(e+x, x)$ is a local $S$-section of $X\times_YX$. Indeed, as $q$ is a $Y$-endomorphism of $X$, the sections $x$ and $q(x)$ (resp. $e+x$ and $q(e+x)=p(e)+q(x)=q(x)$) have the same image in $Y$. 
\smallskip

So $E=0$, as $L$ is basic. So $p$, hence $q$ as well, is an isomorphism.
\smallskip

Let $Z=\mathrm{Ker}(q, \mathrm{Id}_X)$, which is an open and closed sub-$Y$-algebraic space of $X$. If $q(x)=x$, that is, if $x\in Z(S)$, then $Z=X$ and $q=\mathrm{Id}_X$, since $X$ has geometrically irreducible $S$-fibers.

\end{proof}
\smallskip

{\bf Lemma 3.3.} --- \emph{For $i=1,2$, let $X_i, L_i, Y_i=X_i/L_i$ be as in $(2.1)$, with $L_i$ basic, let $x_i\in X_i(S)$ and let $y_i\in Y_i(S)$ be the image of $x_i$. Let $r: Y_1\to Y_2$ be an $S$-isomorphism satisfying $r(y_1)=y_2$.}
\smallskip

\emph{Then there exists a unique $r$-isomorphism $h: X_1\to X_2$ satisfying $h(x_1)=x_2$.}

\begin{proof} Identify $Y_1$ with $Y_2$ by $r$. Write $Y=Y_1=Y_2$ and $y=y_1=y_2$. The uniqueness of $h$ follows from (3.2). Thus, the existence of $h$ is an \'{e}tale local question on $S$. 
\smallskip

One may assume that $S$ is strictly local with closed point $s$. For, if the lemma is proven under this assumption, then $h$ exists in an \'{e}tale neighborhood of $s$ in $S$ by a ``passage \`{a} la limite projective'' (EGA IV 8). As the closed immersion $Y_s\hookrightarrow Y$ induces an equivalence between the category of finite \'{e}tale $Y$-algebraic spaces and the category of finite \'{e}tale $Y_s$-algebraic spaces (SGA 4 XII 5.9 bis), one may further assume that $S=s$ is the spectrum of a separably closed field $k$.
\smallskip

Let $X$ be an algebraic space which is connected and finite \'{e}tale Galois over $Y$ and which dominates $X_1$ and $X_2$, say by $q_i: X\to X_i$, $i=1,2$. By SGA 1 XI 2, $X$ is a (trivial) torsor under a $k$-abelian variety $A$ and $q_i$ is equivariant relative to a unique \'{e}tale $k$-group homomorphism $p_i: A\to A_i$, $i=1,2$. 
\smallskip

With each element $g\in \mathrm{Gal}(X/Y)=G$, there is associated a unique group automorphism $a(g)$ of the $k$-abelian variety $A$ such that $g$ is $a(g)$-equivariant. The map $g\mapsto a(g)$ is a homomorphism of groups and its kernel $a^{-1}(1)=E$ consists of translations by elements of $A(k)$. 
\smallskip

By definition, $\mathrm{Ker}(p_1)$ and $\mathrm{Ker}(p_2)$ are $k$-subgroups of $E_k$. Hence, $X/E$ is dominated by $X_1$ and by $X_2$. So $E_k=\mathrm{Ker}(p_1)=\mathrm{Ker}(p_2)$, as $L_1$ and $L_2$ are basic. 
\smallskip

There exist therefore a unique $Y$-isomorphism $q: X_1\to X_2$ and a unique isomorphism of $k$-abelian varieties $p: A_1\to A_2$ such that $q$ is $p$-equivariant and such that $qq_1=q_2$, $pp_1=p_2$.
\smallskip

Identify $X_1$ with $X_2$ by $q$ and identify $A_1$ with $A_2$ by $p$. Then, $X_1=X_2=X/E$ is Galois over $Y$ with Galois group $G/E$. Clearly, a unique element $h$ of $G/E$ transforms $x_1$ to $x_2$.

\end{proof}
\smallskip

{\bf Proposition 3.4.} --- \emph{Keep the notations of $(2.1)$. Assume that $L$ is basic. Let $f: Y\to S$ denote the structural morphism.}
\smallskip

\emph{Then the $S$-group $G=f_{*}\underline{\mathrm{Aut}}_Y(X)$ is a finite \'{e}tale $S$-algebraic space and the canonical morphism
\[u: G\times_SX\to X\times_YX,\ (g, x)\mapsto (g(x), x)\] is an isomorphism.}

\begin{proof} The finite \'{e}tale $Y$-group $N=\underline{\mathrm{Aut}}_Y(X)$, considered as a proper smooth $S$-algebraic space, has Stein factorization $N\to f_*N\to S$ (SGA 1 X 1.2). By \emph{loc.cit.}, $f_*N\to S$ is finite \'{e}tale and the formation of the Stein factorization commutes with every base change $S'\to S$. 
It remains only to verify that $u\times_Ss$ is an isomorphism for each geometric point $s$ of $S$. One can thus assume that $S$ is the spectrum of a separably closed field. Then, by the proof of (3.3), $X\to Y$ is Galois and hence $u$ is an isomorphism.

\end{proof}
\smallskip

{\bf Proposition 3.5.} --- \emph{For $i=1,2$, let $X_i, L_i, Y_i=X_i/L_i$ be as in $(2.1)$ with $L_i$ basic, let $f_i: Y_i\to S$ be the structural morphism, let $G_i=f_{i*}\underline{\mathrm{Aut}}_{Y_i}(X_i)$ and let $r: Y_1\to Y_2$ be an $S$-isomorphism. Write $X'_i, L'_i, Y'_i, r'$ for the base change of $X_i, L_i, Y_i, r$ by a morphism $S'\to S$. Then the sheaf on $(\mathrm{Sch}/S)$ for the \'{e}tale topology, $S'\mapsto \mathrm{Isom}_{r'}(X'_1, X'_2)$, is a $(G_2, G_1)$-bitorsor. The map
\[\mathrm{Isom}_r(X_1, X_2)\to \mathrm{Isom}_{r'}(X'_1, X'_2),\ h\mapsto h\times_SS'\] is a bijection if $S'\to S$ is $0$-acyclic} (SGA 4 XV 1.3).

\begin{proof} Identify $Y_1$ with $Y_2$ by $r$, write $Y=Y_1=Y_2$ and $f=f_1=f_2$. Let $I=\underline{\mathrm{Isom}}_Y(X_1, X_2)$, which is a finite \'{e}tale $Y$-algebraic space. By SGA 1 X 1.2, $f_*I$ is a finite \'{e}tale $S$-algebraic space and its formation commutes with every base change $S'\to S$ and it is a $(G_2, G_1)$-bitorsor by the same argument of (3.4). When $u: S'\to S$ is $0$-acyclic, the adjunction morphism
\[f_*I\to u_*u^*f_*I\] is an isomorphism and in particular
\[\Gamma(Y, I)=\Gamma(S, f_*I)\to \Gamma(Y', I')=\Gamma(S, u_*u^*f_*I)\] is a bijection.

\end{proof}
\smallskip

4. \emph{Factorization over a separably closed field.}
\smallskip

{\bf Proposition 4.1.} --- \emph{Over a separably closed field $k$, every finite \'{e}tale surjective $k$-morphism from a $k$-abelian variety $A'$ to a $k$-algebraic space $Y$, $A'\to Y$, factors up to $k$-isomorphisms in a unique way as the composite of an \'{e}tale isogeny of $k$-abelian varieties $A'\to A$ and a quotient $A\to A/L=Y$ by a basic finite \'{e}tale $k$-equivalence relation $L$ on $A$.}

\begin{proof} Let $A'\to Y$ be dominated by a finite \'{e}tale Galois $k$-morphism $A''\to Y$ with $A''$ connected. By SGA 1 XI 2, $A''$ may be endowed with a $k$-abelian variety structure so that the projection $A''\to A'$ is an \'{e}tale $k$-group homomorphism. Let $E''$ denote the kernel of this homomorphism. 
\smallskip

Let $E$ consist of the elements of $G=\mathrm{Gal}(A''/Y)$ which are translations of $A''$, namely, of the form $x\mapsto x+a$, with $a\in A''(k)$. 
Evidently, $E$ is a subgroup of $G$, $E\supset E''(k)$ and the quotient $A=A''/E$ inherits a $k$-abelian variety structure from that of $A''$.
\smallskip

The factorization
\[A'=A''/E''\to A=A''/E\to Y=A''/G\] is then as desired. One argues as in (3.3) that such a factorization is unique up to $k$-isomorphisms.

\end{proof}
\smallskip

5. \emph{Basicness is an open and closed property.}
\smallskip

{\bf Proposition 5.1.} --- \emph{Keep the notations of $(2.1)$. Let \[F: (\mathrm{Sch}/S)^o\to (\mathrm{Sets})\] be the following sub-functor of the final functor }:
\smallskip

\emph{For an $S$-scheme $S'$, $F(S')=\{\emptyset\}$, if $L\times_SS'$ is basic, and $F(S')=\emptyset$, otherwise.}
\smallskip

\emph{Then $F$ is representable by an open and closed sub-algebraic space of $S$.}

{\bf Lemma 5.2.} --- \emph{Keep the notations of $(2.1)$. Assume that $S$ is the spectra of a discrete valuation ring with generic point $t$ and closed point $s$.}
\smallskip

\emph{Then the following conditions are equivalent }:
\smallskip

1) \emph{$L_s$ is basic.}
\smallskip

2) \emph{$L_t$ is basic.}

\begin{proof} One may assume that $S$ is strictly henselian.
\smallskip

Assume 1). The finite \'{e}tale morphism $X_s\to X_s/L_s$ is Galois (3.4). So by SGA 4 XII 5.9 bis, $X$ is Galois over $X/L=Y$, as $S$ is strictly henselian and $Y\to S$ is proper. Let $G=\mathrm{Gal}(X/Y)$. Each element $g\in G$ is equivariant relative to a unique group automorphism $a(g)$ of the $S$-abelian algebraic space $A$. The homomorphism $g\mapsto a(g)_s$ is injective, as $L_s$ is basic. It follows that $g\mapsto a(g)_t$ is injective. For, the specialization homomorphism
\[\mathrm{Aut}_{t}(A_t)\ \widetilde{\leftarrow}\ \mathrm{Aut}_S(A)\to \mathrm{Aut}_s(A_s)\] is injective. So $L_t$ is basic.
\smallskip

Assume 2). Let $\overline{t}$ be the spectrum of a separable closure of $k(t)$. By (3.4), $X_{\overline{t}}\to Y_{\overline{t}}$ is Galois with Galois group, say $G$. Replacing if necessary $S$ by its normalization $S'$ in a finite sub-extension of $k(\overline{t})/k(t)$, $X$ by $X\times_SS'$, $Y$ by $Y\times_SS'$ and $L$ by $L\times_SS'$, one may assume that $X_t$ is a $G$-torsor on $Y_t$ and that $L_t$ is defined by the $G$-action on $X_t$. Now $X$ being the $S$-N\'{e}ron model of $X_t$, each $t$-automorphism of $X_t$ uniquely extends to an $S$-automorphism of $X$. In particular, the $G$-action on $X_t$ extends to an action on $X$ with graph evidently equal to that of $L$. For each $g\in G$, let $a(g)$ be the unique $S$-automorphism of the $S$-abelian algebraic space $A$ relative to which $g: X\to X$ is equivariant. The homomorphism $g\mapsto a(g)_t$ is injective, as $L_t$ is basic. Then $g\mapsto a(g)_s$ is injective, for the specialization homomorphism
\[\mathrm{Aut}_t(A_t)\ \widetilde{\leftarrow}\ \mathrm{Aut}_S(A)\to \mathrm{Aut}_s(A_s)\] is injective. So $L_s$ is basic.

\end{proof}
\smallskip

{\bf Lemma 5.3.} --- \emph{Keep the notations of $(2.1)$. Assume that $L$ is basic. Assume that $S$ is an affine scheme and that $(S_i)$ is a projective system of affine noetherian schemes, indexed by a co-directed set $I$, with limit $S$.}
\smallskip

\emph{Then there exist an index $i\in I$, an $S_i$-abelian algebraic space $A_i$, an $A_i$-torsor $X_i$ on $S_i$ for the \'{e}tale topology and a basic finite \'{e}tale $S_i$-equivalence relation $L_i$ on $X_i$ such that $A=A_i\times_{S_i}S$, $X=X_i\times_{S_i}S$ and $L=L_i\times_{S_i}S$.}

\begin{proof} By the technique of ``passage \`{a} la limite projective'' (EGA IV 8, 9), there exist $i_o\in I$ and $A_{i_o}, X_{i_o}, L_{i_o}$ as desired except possibly the property of being \emph{basic}. By (5.2), for each $j\geq i_o$, $L_{i_o}\times_{S_{i_o}}S_j$ is basic at precisely the points of an open and closed sub-scheme $S'_j$ of $S_j$. The projection $S\to S_j$ factors through $S'_j$, as $L$ is basic. Now $(S_j-S'_j)_{j\geq i_o}$ is a projective system of affine noetherian schemes with empty limit. So $S_i=S'_i$ for some $i\geq i_o$.

\end{proof}
\smallskip

5.4. \emph{Proof of} (5.1).
\smallskip

The question being an \'{e}tale local question on $S$, one may assume that $S$ is a scheme, then affine and by (5.3) noetherian. The claim is now immediate by (5.2).

\begin{flushright}
$\square$
\end{flushright}

\smallskip

6. \emph{Automorphisms II.}
\smallskip

We study the $S$-automorphisms of $Y$.
\smallskip

{\bf Lemma 6.1.} --- \emph{Keep the notations of $(2.1)$. Assume that $L$ is basic. Let $p: X\to X/L=Y$ denote the projection. Let $r$ be an $S$-automorphism of $Y$.}
\smallskip

\emph{Then there is at most one $S$-section $a$ of $A$ which verifies $rp=pT_a$, where $T_a: x\mapsto a+x$ is the $S$-automorphism of $X$, the translation by $a$.}

\begin{proof} If $rp=pT_a=pT_b$ holds for two $S$-sections $a, b$ of $A$, then $T_{a-b}$ acts on $X$ as a $Y$-automorphism. In particular, $T_{a-b}$ is locally on $S$ of finite order. So $T_{a-b}=\mathrm{Id}_X$ and $a=b$, as $L$ is basic.

\end{proof}

\smallskip

{\bf Lemma 6.2.} --- \emph{Keep the notations and assumptions of $(6.1)$. }
\smallskip

\emph{Let $F: (\mathrm{Sch}/S)^o\to (\mathrm{Sets})$ be the following sub-functor of the final functor }:
\smallskip

\emph{For an $S$-scheme $S'$, $F(S')=\{\emptyset\}$, if $r_{S'}p_{S'}=p_{S'}T_{a'}$ for an $S'$-section $a'$ of $A_{S'}$, and $F(S')=\emptyset$, otherwise.}
\smallskip

\emph{Then $F$ is representable by an open and closed sub-algebraic space of $S$.}

\begin{proof} The question being by (6.1) an \'{e}tale local question on $S$, one may assume $S$ to be a scheme, then affine and by (5.3) noetherian.
\smallskip

--- \emph{The functor $F$ verifies the valuative criterion of properness }:
\smallskip

Namely, given an $S$-scheme $S'$ which is the spectra of a discrete valuation ring and which has generic point $t'$, then $F(S')=F(t')$. For, if $r_{t'}p_{t'}=p_{t'}T_{a_{t'}}$ holds for some point $a_{t'}\in A(t')$, then $a_{t'}$ extends uniquely to an $S'$-section $a'$ of $A_{S'}$ and the equation $r_{S'}p_{S'}=p_{S'}T_{a'}$ holds.
\smallskip

--- \emph{The functor $F$ is formally \'{e}tale }:
\smallskip

Namely, $F(S')=F(S'')$ for every nilpotent $S$-immersion $S''\hookrightarrow S'$. Indeed, assume $r_{S''}p_{S''}=p_{S''}T_{a''}$ holds for some section $a''\in A(S'')$. As both $p$ and $rp$ are \'{e}tale, there is a unique $S'$-automorphism $T'$ of $X_{S'}$ such that $T'$ restricts to $T_{a''}$ on $X_{S''}$ and such that $r_{S'}p_{S'}=p_{S'}T'$ holds. This $T'$ is equivariant with respect to a unique $S'$-group automorphism $\varphi'$ of $A_{S'}$. As the $S$-group $\underline{\mathrm{Aut}}_S(A)$ is unramified over $S$ and as $\varphi'$ restricts to the identity automorphism of $A_{S''}$, it follows that $\varphi'=\mathrm{Id}_{A_{S'}}$. So $T'$ is of the form $T_{a'}$ for an $S'$-section $a'\in A(S')$.
\smallskip

It is now evident that $F\to S$ is representable by an open and closed immersion.

\end{proof}

\smallskip

{\bf Proposition 6.3.} --- \emph{Keep the notations of $(2.1)$. Assume that $L$ is basic. Let $p: X\to X/L=Y$ denote the projection and $f: Y\to S$ the structural morphism.}
\smallskip

\emph{Then }:
\smallskip

\emph{a) The $S$-subgroup $R$ of $\underline{\mathrm{Aut}}_S(Y)$ defined to consist of the local $S$-automorphisms $r$ satisfying $rp=pT_a$ for some local $S$-sections $a$ of $A$ is open and closed in $\underline{\mathrm{Aut}}_S(Y)$.}
\smallskip

\emph{b) The homomorphism $r\mapsto a$ is a closed immersion of $R$ into $A$. In particular, $R$ is commutative and representable by a proper $S$-algebraic space.}
\smallskip

\emph{c) An $S$-section $a$ of $A$ lies in the image of $R$ if and only if $T_ag=gT_a$ for all local $S$-sections $g$ of $G=f_*\underline{\mathrm{Aut}}_Y(X)$.}

\begin{proof} Note that $R\to A$, $r\mapsto a$, is a well-defined morphism by (6.1). This morphism is a monomorphism because $p$, being \'{e}tale surjective, is an epimorphism in the category of $S$-algebraic spaces.
\smallskip

The claim that the sub-$S$-group functor $R\subset\underline{\mathrm{Aut}}_S(Y)$ is open and closed is a rephrase of (6.2), hence \emph{a}).
\smallskip

If an $S$-section $a$ of $A$ satisfies $T_ag=gT_a$ for every local $S$-section $g$ of $G$, then $T_a$, being $G$-equivariant, defines by passing to the $G$-quotient an $S$-automorphism $r$ of $X/G=Y$ (3.4). One has thus $rp=pT_a$ by construction. Conversely, the identity $rp=pT_a$ implies that, for every local $S$-section $g$ of $G$, one has $rp=rpg^{-1}=pgT_ag^{-1}$ and then by (6.1) $gT_ag^{-1}=T_a$. The characterization in \emph{c}) therefore follows and it implies evidently that the monomorphism $R\hookrightarrow A$, $r\mapsto a$, is a closed immersion, hence \emph{b}).

\end{proof}

\smallskip

7. \emph{Rigidity of torsors and albanese.}
\smallskip

{\bf Lemma 7.1.} --- \emph{Let $S$ be the spectra of a discrete valuation ring and $t$ (resp. $s$) a geometric generic (resp. closed) point of $S$. Let $X$ be a proper smooth $S$-algebraic space.}
\smallskip

\emph{Then the following conditions are equivalent }:
\smallskip

1) \emph{$X\times_Ss$ has an $s$-abelian variety structure.}
\smallskip

2) \emph{$X\times_St$ has a $t$-abelian variety structure.}

\begin{proof} One may assume that $S$ is strictly henselian. Then $X$ has $S$-sections by EGA IV 17.16.3. Fix an $S$-section $o$ of $X$. When 1) (resp. 2)) holds, after a translation of the origin, $X\times_Ss$ (resp. $X\times_St$) has an abelian variety structure over $s$ (resp. $t$) with $o_s$ (resp. $o_t$) as its zero section. Under either assumption, $X$ has geometrically irreducible $S$-fibers and admits a unique $S$-abelian algebraic space structure with $o$ being the zero section (``Geometric Invariant Theory'' 6.14).

\end{proof}

\smallskip

{\bf Proposition 7.2.} --- \emph{Let $S$ be an algebraic space and $X$ a proper smooth $S$-algebraic space. Let $T: (\mathrm{Sch}/S)^o\to (\mathrm{Sets})$ be the following sub-functor of the final functor }:
\smallskip

\emph{For an $S$-scheme $S'$, $T(S')=\{\emptyset\}$, if $X\times_Ss'$ has an $s'$-abelian variety structure for each geometric point $s'$ of $S'$, and $T(S')=\emptyset$, otherwise.}
\smallskip

\emph{Then $T$ is representable by an open and closed sub-algebraic space of $S$. And $X\times_ST=X_T$ is a $\mathrm{Pic}^o_{\mathrm{P}_T/T}$-torsor on $T$ for the \'{e}tale topology where $\mathrm{P}_T=\mathrm{Pic}^o_{X_T/T}$ is a $T$-abelian algebraic space.}

\begin{proof} From (7.1) and by a ``passage \`{a} la limite'', one deduces that $T\to S$ is representable by an open and closed immersion.
Replacing $S$ by $T$, one may assume $T=S$. Let $S'\to S$ be a smooth morphism with geometrically irreducible fibers such that $X'=X\times_SS'$ admits an $S'$-section $e'$. One may for instance take $S'=X$ and $e'$ to be the diagonal section $\Delta_{X/S}$. Notice that in $S''=S'\times_SS'$ the only open and closed neighborhood of the diagonal is $S''$.
By ``Geometric Invariant Theory'' 6.14, on $X'$ there is a unique $S'$-abelian algebraic space structure with zero section $e'$. It follows that $\mathrm{P}=\mathrm{Pic}^{\tau}_{X/S}$, being representable by an $S$-algebraic space (Artin), is an $S$-abelian algebraic space and that $\mathrm{Pic}^{\tau}_{X/S}=\mathrm{Pic}^o_{X/S}$, as one verifies after the \emph{fppf} base change $S'\to S$ on the $S'$-abelian algebraic space $X'$. Let the dual $S$-abelian algebraic space of $\mathrm{P}$ be $A$. Let $p_1, p_2$ be the two projections of $S''$ onto $S'$ and $u: p_1^*X'\to p_2^*X'$ the descent datum on $X'$ relative to $S'\to S$ corresponding to $X$. By (3.1)+(7.3) $u$ is equivariant with respect to a unique $S''$-group automorphism $a$ of $A''=A\times_SS''$. It suffices evidently to show that $a=1$. Now the $S$-group $\underline{\mathrm{Aut}}_S(A)$ being unramified and separated over $S$, the relation ``$a=1$'' is an open and closed relation on $S''$ and holds on the diagonal and so holds.

\end{proof}

\smallskip

{\bf Lemma 7.3.} --- \emph{Let $S$ be an algebraic space, $A$ an $S$-abelian algebraic space and $X$ an $A$-torsor on $S$ for the \'{e}tale topology.}
\smallskip

\emph{Then there exists a canonical isomorphism
\[X\stackrel{A}{\wedge} \mathrm{Pic}_{A/S}\ \widetilde{\to}\ \mathrm{Pic}_{X/S}\] which induces isomorphisms
\[\mathrm{Pic}^o_{A/S}=X\stackrel{A}{\wedge}\mathrm{Pic}^o_{A/S}\ \widetilde{\to}\ \mathrm{Pic}^o_{X/S},\]
\[\mathrm{NS}_{A/S}=X\stackrel{A}{\wedge}\mathrm{NS}_{A/S}\ \widetilde{\to}\ \mathrm{NS}_{X/S}.\]}

\begin{proof} This is ``Faisceaux amples sur les sch\'{e}mas en groupes et les espaces homog\`{e}nes'' XIII 1.1, by which it is justified to call $A$ the \emph{albanese} of $X$.

\end{proof}

\smallskip

8. \emph{Specialization and descent}.

\smallskip

{\bf Proposition 8.1.} --- \emph{Let $S$ be the spectra of a discrete valuation ring and $t$ (resp. $s$) a geometric generic (resp. closed) point of $S$. Let $Y$ be a proper smooth $S$-algebraic space. Suppose that $Y_{t}$ is the quotient of a $t$-abelian variety $A_t$ by a basic finite \'{e}tale $t$-equivalence relation such that $[A_t: Y_t]$ is prime to the characteristic of $k(s)$.}
\smallskip

\emph{Then $Y_s$ is the quotient of an $s$-abelian variety by a basic finite \'{e}tale $s$-equivalence relation.}

\begin{proof} One may assume that $S$ is strictly henselian with closed point $s$. Recall that by (3.4) $A_t$ is Galois over $Y_t$. As $Y$ is proper smooth over $S$, $S$ strictly henselian and $G=\mathrm{Gal}(A_t/Y_t)$ of order prime to the characteristic of $k(s)$, the specialization morphism
\[H^1(Y, G)\to H^1(Y_t, G)\] is by SGA 4 XVI 2.2 a bijection. One finds thus a $G$-torsor on $Y$ for the \'{e}tale topology, $A\to Y$, which specializes to $A_t\to Y_t$ at $t$. In particular, $A$ is proper smooth over $S$ and has $S$-sections by EGA IV 17.16.3 and has by (7.2) an $S$-abelian algebraic space structure. The finite \'{e}tale $S$-equivalence relation on $A$ of graph $A\times_YA$ is basic because it is basic at $t$ (5.1).

\end{proof}

\smallskip

{\bf Proposition 8.2.} --- \emph{Let $S$ be an algebraic space. Let $Y$ be an algebraic space which is proper flat of finite presentation over $S$.}
\smallskip

\emph{Let $U: (\mathrm{Sch}/S)^o\to (\mathrm{Sets})$ be the following sub-functor of the final functor }:
\smallskip

\emph{For an $S$-scheme $S'$, $U(S')=\{\emptyset\}$, if $Y\times_Ss'$ is the quotient of $s'$-abelian variety by a finite \'{e}tale $s'$-equivalence relation for each geometric point $s'$ of $S'$, and $U(S')=\emptyset$, otherwise.}
\smallskip

\emph{Then $U\to S$ is representable by an open immersion of finite presentation. There exist a $U$-abelian algebraic space $A$, an $A$-torsor $X$ on $U$ for the \'{e}tale topology and a basic finite \'{e}tale $U$-equivalence relation $L$ on $X$ such that $Y\times_SU$ is $U$-isomorphic to the quotient $X/L$.}

\begin{proof} Observe that, by (4.1), for a given geometric point $s$ of $S$, $U(s)=\{\emptyset\}$ if and only if $Y_s$ is the quotient of an $s$-abelian variety by a \emph{basic} finite \'{e}tale $s$-equivalence relation.
\smallskip

--- \emph{Reduction to the case where $Y$ is $S$-smooth }:
\smallskip

Consider the functor $V: (\mathrm{Sch}/S)^o\to (\mathrm{Sets})$ :
\smallskip

\emph{For an $S$-scheme $S'$, $V(S')=\{\emptyset\}$, if $Y\times_Ss'$ is smooth over $s'$ for each geometric point $s'$ of $S'$, and $V(S')=\emptyset$, otherwise.}
\smallskip

It is evident that $V\to S$ is representable by an open immersion of finite presentation, that $Y\times_SV$ is proper and smooth over $V$ and that $V$ contains $U$ as a sub-functor.
Restricting to $V$, one may assume that $S=V$, namely, that $Y$ is smooth over $S$. 
\smallskip

--- \emph{Assume that $Y$ is $S$-smooth. Then $U\to S$ is representable by an open immersion of finite presentation. There exist an \'{e}tale surjective morphism $U'\to U$, a $U'$-abelian algebraic space $A'$, an $A'$-torsor $X'$ on $U'$ for the \'{e}tale topology and a basic finite \'{e}tale $U'$-equivalence relation $L'$ on $X'$ such that $Y\times_SU'$ is $U'$-isomorphic to the quotient $X'/L'$ }:
\smallskip

The question being an \'{e}tale local question on $S$, one may assume that $S$ is a scheme, then affine, then by (5.3) noetherian and strictly local with closed point $s$ and finally that $U(s)=\{\emptyset\}$, i.e., that $Y_s$ is the quotient of an $s$-abelian variety $A_s$ by a basic finite \'{e}tale $s$-equivalence relation $L_s$. By SGA 4 XII 5.9 bis, there is a finite \'{e}tale $S$-morphism $A\to Y$ which specializes to $A_s\to A_s/L_s=Y_s$ at $s$. This algebraic space $A$ is in particular proper smooth over $S$ and has closed fiber the abelian variety $A_s$. As $S$ is strictly local, $A$ admits $S$-sections (EGA IV 17.16.3). So by (7.2) $A$ has an $S$-abelian algebraic space structure. The finite \'{e}tale $S$-equivalence relation $L$ on $A$ of graph $A\times_YA$, basic at $s$, is basic (5.1).
\smallskip

--- \emph{Assume that $Y$ is $S$-smooth. Then up to an \'{e}tale localization on $U'$ there exists a descent datum on $(A', X', L')$ relative to $U'\to U$ }:
\smallskip

Restricting to the open sub-algebraic space $U$, one may assume that $U=S$. Write $S'$ for $U'$ and $Y'$ for $Y\times_SS'$. It suffices to prove the existence of a finite \'{e}tale $Y$-algebraic space $X$ such that the finite \'{e}tale $Y'$-algebraic space 
$\underline{\mathrm{Isom}}_{Y'}(X\times_YY', X')$ is surjective over $Y'$, for then $X$ is an $A$-torsor on $S$ for the \'{e}tale topology (7.2) where the $S$-abelian algebraic space $A$ satisfies $\mathrm{Pic}^o_{A/S}=\mathrm{Pic}^o_{X/S}$ (7.3) and $X\times_YX$ is the graph of a basic finite \'{e}tale $S$-equivalence relation $L$ on $X$.
\smallskip

i) \emph{Case where $Y$ has an $S$-section $y$ }:
\smallskip

By an \'{e}tale localization on $S'$, one may assume that there is an $S'$-section $x'$ of $X'$ which is mapped to $y'=y\times_SS'$ by $X'\to Y'$. Let $S''=S'\times_SS'$ and $p_1, p_2$ the two projections of $S''$ onto $S'$. By (3.3), there exists a unique $(Y\times_SS'')$-isomorphism $h: p_1^*X'\to p_2^*X'$ which transforms $p_1^*(x')$ to $p_2^*(x')$. That is, $h$ is a gluing datum on $(X', x')$ relative to $S'\to S$. By (3.3) again, $h$ is a descent datum. This descent provides a finite \'{e}tale $Y$-algebraic space $X$ which verifies $X\times_YY'=X'$ and which is equipped with an $S$-section $x$ having image $y$ in $Y(S)$.
\smallskip

ii) \emph{General case }:
\smallskip

By the second projection, $Y\times_SY=Y_1$ has a $Y$-algebraic space structure which admits the $Y$-section $\Delta_{Y/S}$. One finds by i) a finite \'{e}tale $Y_1$-algebraic space $X_1$ such that the finite \'{e}tale $Y'_1$-algebraic space
\[\underline{\mathrm{Isom}}_{Y'_1}(X_1\times_{Y_1}Y'_1, X'\times_{Y'}Y'_1)\] is surjective over $Y'_1$, where $Y'_1:=Y\times_SY'=Y_1\times_YY'$. Let the Stein factorization of the proper smooth $Y$-algebraic space $X_1$ be 
\[X_1\to X\to Y\] and let $c$ be the canonical morphism
\[c: X_1\to X\times_YY_1.\] By SGA 1 X 1.2, the morphism $X\to Y$ is finite \'{e}tale and the formation of the Stein factorization commutes with every base change $T\to S$. It suffices to show that $c$ is an isomorphism, for then 
\[\underline{\mathrm{Isom}}_{Y'_1}(X_1\times_{Y_1}Y'_1, X'\times_{Y'}Y'_1)=\underline{\mathrm{Isom}}_{Y'}(X\times_YY', X')\times_{Y'}Y'_1\] and $X$ is the sought after $Y$-algebraic space. This amounts to showing that $c\times_Ss$ is an isomorphism for each geometric point $s$ of $S$. One may thus assume that $S$ is the spectrum of a separably closed field. Then $S'$, being \'{e}tale surjective over $S$, is a non-empty disjoint union of $S$. Index these components of $S'$ as $S_i$, $i\in\pi_o(S')=\Pi$, write $X'=\sum_{i\in\Pi}X_i$, fix a point $y\in Y(S)$ and choose a point $x_i\in X_i(S_i)$ above $y$ for each $i\in\Pi$. These $(X_i, x_i)$'s are all mutually $Y$-isomorphic by (3.3). Clearly, $c$ is an isomorphism.

\end{proof}

\smallskip

9. \emph{Proof of Theorem} 1.1.
\smallskip

\emph{a}) Let the open sub-algebraic space $U$ of $S$, the $U$-abelian algebraic space $A$, the $A$-torsor $X$ on $U$ for the \'{e}tale topology be as in (8.2) so that $Y\times_SU=Y_U$ is the quotient of $X$ by a basic finite \'{e}tale $U$-equivalence relation $L$. It suffices to show that $U=S$. 
\smallskip

Factor $A'_t\to Y_t$ as $A'_t\to X_t\to Y_t$ (4.1). The degree $[X_t :Y_t]=d$, which divides $[A'_t: Y_t]$, is prime to the residue characteristics of $S$. By (3.4), there is a maximal open and closed sub-algebraic space $U'$ of $U$ such that $X_{U'}$ is of constant degree $d$ over $Y_{U'}$. By (8.1), $U'$ is closed in $S$. So $U'=U=S$, as $S$ is connected. 
\smallskip

\emph{b}) The assertion is immediate by (3.5)+(7.3). And the map which with $P'\in H^1(S, G)$ associates 
\[X'=P'\stackrel{G}{\wedge}X\] provides the inverse.

\begin{flushright}
$\square$
\end{flushright}

\smallskip

10. \emph{Irreducibility of non-uniruled degenerate fibers. Almost non-degeneration}.
\smallskip

The irreducibility of non-uniruled degenerate fibers of an abelian fibration is shown in \cite{almost}, hence :
\smallskip

{\bf Proposition 10.1.} --- \emph{Keep the notations of $(1.1)$. Assume that $S$ is open in an algebraic space $\overline{S}$ and that $Y$ is open dense in a separated finitely presented $\overline{S}$-algebraic space $\overline{Y}$ with structural morphism $\overline{f}$. Let $\overline{s}$ be a geometric point of $\overline{S}$ with values in an algebraically closed field.}
\smallskip

\emph{Then $\overline{f}^{-1}(\overline{s})$ is irreducible if it is non-empty and does not have uniruled irreducible components and if one of the following two conditions holds }:
\smallskip

\emph{a) $\overline{f}$ is flat at all maximal points of $\overline{f}^{-1}(\overline{s})$.}
\smallskip

\emph{b) $\overline{S}$ is locally noetherian and $\overline{f}$ is universally open at all maximal points of $\overline{f}^{-1}(\overline{s})$.}

\begin{proof} By standard arguments one may assume that $\overline{S}$ is the spectra of a complete discrete valuation ring with $\overline{S}-S=\{\overline{s}\}$, that $\overline{f}$ is flat and that $\overline{f}^{-1}(\overline{s})$ has no imbedded components. Let $X$ be as in (1.1) and $\overline{X}$ the normalization of $\overline{Y}$ in $X$. Then $\overline{X}_{\overline{s}}$ is non-empty and does not have uniruled irreducible components and hence by \cite{almost} 4.1 is irreducible. So $\overline{f}^{-1}(\overline{s})$ is irreducible.

\end{proof}

\smallskip

{\bf Theorem 10.2.} --- \emph{Let $S$ be the spectra of a discrete valuation ring and $\overline{t}$ (resp. $\overline{s})$ a geometric generic (resp. closed) point of $S$ with values in an algebraically closed field. Let $Y$ be an $S$-algebraic space with structural morphism $f$. Assume that $f$ is separated of finite type, that $Y$ is normal integral and at each of its geometric codimension $\geq 2$ points either regular or pure geometrically para-factorial of equal characteristic, that $f^{-1}(\overline{t})$ is the quotient of a $\overline{t}$-abelian variety $A'_{\overline{t}}$ by a finite \'{e}tale $\overline{t}$-equivalence relation with degree $[A'_{\overline{t}}: f^{-1}(\overline{t})]$ prime to the characteristic of $k(\overline{s})$ and that $f^{-1}(\overline{s})$ is non-empty, proper, of total multiplicity prime to the characteristic of $k(\overline{s})$ and does not have uniruled irreducible components.}
\smallskip

\emph{Then $f$ factors canonically as $Y\to E\to S$ with $Y\to E$ proper and smooth, where $E$ is a finite flat $S$-algebraic stack, regular, tame over $S$ and satisfies $E\times_S\overline{t}=\overline{t}$.}

\begin{proof} Notice that $f$ is faithfully flat and, being separated of finite type with geometrically irreducible (10.1) and proper fibers, that $f$ is also proper (EGA IV 15.7.10). Let the total multiplicity of $f^{-1}(\overline{s})$, that is by definition, the greatest common divisor of the lengths of the local rings of $f^{-1}(\overline{s})$ at its maximal points, be $\delta$, which by hypothesis is prime to the characteristic of $k(\overline{s})$. Thus, the $1$-codimensional cycle on $Y$ with rational coefficients, $\Delta=f^*\mathrm{Cyc}_S(\pi)/\delta$, where $\pi$ is a uniformizer of $S$, has integral coefficients and is a prime cycle and is locally principal, for $Y$, being normal and geometrically para-factorial at all its geometric codimension $\geq 2$ points, has geometrically factorial local rings (EGA IV 21.13.11). With $\Delta$ one associates a canonical $\mu_{\delta}$-torsor on $Y$ for the \'{e}tale topology, $Y'\to Y$. Let $S'=\mathrm{Spec}\ \Gamma(Y', \mathcal{O}_{Y'})$. There is by quotient by $\mu_{\delta}$ an $S$-morphism
\[Y=[Y'/\mu_{\delta}]\to [S'/\mu_{\delta}]=E.\] 

It suffices evidently to show that $Y'\to S'$ is smooth, for then 
\[Y\to E\to S\] is the desired factorization of $f$. 
Replacing $f$ by $Y'\to S'$, we assume from now on that $\delta=1$, namely, that $f^{-1}(\overline{s})$ is integral. And, replacing $f$ by $f\times_SS_{(\overline{s})}$, where $S_{(\overline{s})}$ is the strict henselization of $S$ at $\overline{s}$, we assume that $S$ is strictly henselian.
\smallskip

Let $t$ (resp. $s$) be the generic (resp. closed) point of $S$. Choose by (1.1) a triple $(A_t, X_t, L_t)$ so that, for a $t$-abelian variety $A_t$, $Y_t=f^{-1}(t)$ is the quotient of an $A_t$-torsor $X_t$ by a basic finite \'{e}tale $t$-equivalence relation $L_t$. Let $G_t=f_{t*}\underline{\mathrm{Aut}}_{Y_t}(X_t)$. The degree $[X_t: f^{-1}(t)]=[G_t: t]$, which divides $[A'_{\overline{t}}: f^{-1}(\overline{t})]$ (4.1), is prime to the characteristic of $k(s)$.
\smallskip

Let $X$ be the normalization of $Y$ in $X_t$. Then $X_{\overline{s}}$ is non-empty and does not have uniruled irreducible components and hence is irreducible (10.1). Let $p: X\to Y$ be the projection and $x$ (resp. $y=p(x)$) the generic point of $X_s$ (resp. $Y_s$). Note that there is an open neighborhood of $x$ (resp. $y$) in $X$ (resp. $Y$) which is a scheme (\cite{raynaud_specialization} 3.3.2). 
\smallskip

--- \emph{Reduction to the case where $G_t$ is constant and cyclic }:
\smallskip

One has that
$[X_t: Y_t]=[\mathcal{O}_{X, x}: \mathcal{O}_{Y, y}]=e[k(x): k(y)]$ is prime to the characteristic of $k(s)$, where $e$ is the ramification index of $\mathcal{O}_{X, x}$ over $\mathcal{O}_{Y, y}$. One deduces that there is an open sub-scheme $V$ of $Y$ containing $y$ such that $V$ is $S$-smooth, that $(U_s)_{\mathrm{red}}$ is finite \'{e}tale surjective over $V_s$ of rank $[k(x): k(y)]$, where $U=p^{-1}(V)$, and that the ideal of $U$ defining the closed sub-scheme $(U_s)_{\mathrm{red}}$ is generated by one section $h\in\Gamma(U, \mathcal{O}_U)$. In particular, $U$ is regular with $h_u$ being part of a regular system of parameters at each point $u$ of $U_s$. Now $S$ being strictly henselian, choose a point $u\in U_s(s)$, let $v=p(u)$, let $n=\mathrm{dim}_u(U_s)$ and choose $h_1,\cdots, h_n\in\mathcal{O}_{U, u}$ so that $\{h_1\ \mathrm{mod}\ h,\cdots, h_n\ \mathrm{mod}\ h\}$ is the image of a regular system of parameters of $V_s$ at $v$. Then $h, h_1,\cdots, h_n$ form a regular system of parameters of $U$ at $u$. Let 
\[R=\mathrm{Spec}(\mathcal{O}_{U, u}/(h_1,\cdots, h_n)),\] which is regular local of dimension $1$ and finite flat tame along $s$ of rank $e$ over $S$. Let $r$ be the generic point of $R$. The closed image of $p(r)$ in $Y$ is an $S$-section lying in $V$ and one obtains the following commutative diagram of $S$-schemes :
{\[\xymatrix{
                  R \ar[r]^{} \ar[d]_{}  & U=p^{-1}(V) \ar[d]^{p} \\
                  S \ar[r]_{}                & V}\] }The $G_t$-torsor structure on $p^{-1}(p(r))$, where one identifies the $t$-rational point $p(r)$ with $t$, induces an epimorphism 
\[G_t\times_tr\to p^{-1}(p(r)),\ (g, \lambda)\mapsto g.\lambda\] and an isomorphism 
\[G_t/\underline{\mathrm{Norm}}_{G_t}(r)\ \widetilde{\to}\ \pi_o(p^{-1}(p(r))).\] Let $Z_t=X_t/N$ be the quotient of $X_t$ by $N=\underline{\mathrm{Norm}}_{G_t}(r)$ and let $Z$ be the normalization of $Y$ in $Z_t$. The fiber $Z_s$ is irreducible with generic point $z$ being the image of $x$. Observe that, on writing $w$ for the image of $u$ in $Z$, $Z$ is by construction \'{e}tale over $Y$ at $w$ and \emph{a priori} at $z$. So $Z\to Y$ is \'{e}tale, as $Y$ by hypothesis is pure at all its geometric codimension $\geq 2$ points. Observe next that $N$ is constant and cyclic. Replacing $Y$ by $Z$ and $G_t$ by $N$, one may assume that $G_t\simeq \mathbf{Z}/e\mathbf{Z}$.
\smallskip

--- \emph{Assume $G_t=\mathbf{Z}/e\mathbf{Z}$. Reduction to the case where $p: X\to Y$ is \'{e}tale }:
\smallskip

Identify $\mathbf{Z}/e\mathbf{Z}$ with $\mu_e$. The $\mu_e$-torsor $X_t\to Y_t$ corresponds to a pair $(J_t, \alpha_t)$ which consists of an invertible module $J_t$ on $Y_t$ and of an isomorphism $\alpha_t: \mathcal{O}_{Y_t}\ \widetilde{\to}\ J_t^{\otimes e}$. There is an invertible $\mathcal{O}_Y$-module $J$ extending $J_t$, since $Y$ has geometrically factorial local rings. Since furthermore $f$ has geometrically integral fibers, $J^{\otimes e}$ is isomorphic to $\mathcal{O}_Y$, say by $\beta: \mathcal{O}_Y\ \widetilde{\to}\ J^{\otimes e}$. The difference in $H^1(Y_t, \mu_e)$ of the classes of $(J_t, \alpha_t)$ and $(J_t, \beta_t)$, that is, the class of $(\mathcal{O}_{Y_t}, \beta_t^{-1}\alpha_t)$, is contained in the image of the map
\[f^*: H^1(t, \mu_e)\to H^1(Y_t, \mu_e),\] for one has $\Gamma(Y_t, \mathcal{O}_{Y_t})=k(t)$. By (1.1) \emph{b}) replacing $X_t$ by the $\mu_e$-torsor $X'_t\to Y_t$ defined by $(J_t, \beta_t)$, one may assume that $p: X\to Y$ is \'{e}tale. It suffices to show that $X$ is smooth over $S$. This follows from \cite{almost} 4.3.

\end{proof}

\smallskip

{\bf Lemma 10.3.} --- \emph{Keep the notations of $(1.1)$. Assume that $S$ is a noetherian local scheme with closed point $s$ such that one of the following two conditions holds }:
\smallskip

\emph{a) $S$ is regular of dimension $>0$.}
\smallskip

\emph{b) $S$ is pure geometrically para-factorial along $s$ of equal characteristic.}
\smallskip

\emph{Let $U=S-\{s\}$. Then each $U$-section of $f|U$ extends uniquely to an $S$-section of $f$.}

\begin{proof} In case \emph{a}) one applies \cite{almost} 2.1 as the geometric fibers of $f$ do not contain rational curves. Assume \emph{b}). One may assume $S$ strictly local. Choose $(A, X, L)$ as in (1.1), let $G=f_{*}\underline{\mathrm{Aut}}_Y(X)$, $p: X\to Y$ the projection and $y: U\to Y$ a section of $f|U$. The finite \'{e}tale $S$-group $G$ is constant and the $G|U$-torsor $p^{-1}(y)$ is trivial, for $S$ is strictly local and pure along $s$. By \emph{loc.cit.} 5.2+5.3 each $U$-section of $p^{-1}(y)$ extends uniquely to an $S$-section of $X$, hence the claim.

\end{proof}

\smallskip

{\bf Lemma 10.4.} --- \emph{Let $S$ be a noetherian normal local scheme of equal characteristic zero pure geometrically para-factorial along its closed point $s$. Let $U=S-\{s\}$. Let $E$ be the fiber category on the category of $S$-algebraic spaces whose fiber over each $S$-algebraic space $S'$ is the full sub-category of the category of $S'$-algebraic spaces consisting of the $S'$-algebraic spaces $Y'$ which, for an $S'$-abelian algebraic space $A'$, an $A'$-torsor $X'$ on $S'$ for the \'{e}tale topology of finite order and a basic finite \'{e}tale $S'$-equivalence relation $L'$ on $X'$, is representable as the quotient $X'/L'$.}

\smallskip

\emph{Then the restriction functor $E(S)\to E(U)$ is an equivalence of categories.}

\begin{proof} This restriction functor is fully faithful by (10.3), since every $S$-smooth algebraic space is pure geometrically para-factorial along its closed $S$-fiber. 
\smallskip

This functor is essentially surjective. For, $U$-abelian algebraic spaces extend to $S$-abelian algebraic spaces (\cite{almost} 5.1). And, as $S$ is pure along $s$, if $n$ is an integer $\geq 1$, $A$ an $S$-abelian algebraic space and ${}_nA=\mathrm{Ker}(n.\mathrm{Id}_A)$, ${}_nA|U$-torsors on $U$ for the \'{e}tale topology extend to ${}_nA$-torsors on $S$ for the \'{e}tale topology. Thus, every object of $E(U)$ is a quotient $(X|U)/(G|U)$ where, for an $S$-abelian algebraic space $A$, $X$ is an $A$-torsor on $S$ for the \'{e}tale topology of finite order and $G$ is a finite \'{e}tale $S$-group such that $G|U$ acts and defines a basic finite \'{e}tale $U$-equivalence relation on $X|U$. Each such action $G|U\times_UX|U\to X|U$ extends by \cite{almost} 5.2+5.3 to a unique $S$-morphism $\mu: G\times_SX\to X$, since $G\times_SX$ is geometrically para-factorial along its closed $S$-fiber. Clearly, $\mu$ represents a $G$-action and defines a basic finite \'{e}tale $S$-equivalence relation on $X$. And $X/G$ is the desired extension of $(X|U)/(G|U)$.

\end{proof}

\smallskip

{\bf Theorem 10.5.} --- \emph{Let $S$ be an integral scheme with generic point $t$. Let $F$ be the fiber category on the category of $S$-algebraic spaces whose fiber over each $S$-algebraic space $S'$ is the full sub-category of the category of $S'$-algebraic spaces consisting of the $S'$-algebraic spaces $Y'$ which, for an $S'$-abelian algebraic space $A'$, an $A'$-torsor $X'$ on $S'$ for the \'{e}tale topology and a basic finite \'{e}tale $S'$-equivalence relation $L'$ on $X'$, is representable as the quotient $X'/L'$. Let $Y$ be an $S$-algebraic space with structural morphism $f$. Assume that $Y$ is locally noetherian normal integral of residue characteristics zero and at all its geometric codimension $\geq 2$ points pure and geometrically para-factorial. Assume furthermore that $f^{-1}(t)$ is an object of $F(t)$ and that, for each geometric codimension $1$ point $\overline{y}$ of $Y$, $f\times_SS_{(\overline{s})}$ is separated of finite type and flat at $\overline{y}$ and the geometric fiber $f^{-1}(\overline{s})$ is proper and does not have uniruled irreducible components, where $S_{(\overline{s})}$ denotes the strict localization of $S$ at $\overline{s}=f(\overline{y})$.}
\smallskip

\emph{Then up to unique isomorphisms there exists a unique groupoid in the category of $S$-algebraic spaces whose nerve $(Y., d_., s_.)$ satisfies the following conditions }:
\smallskip

\emph{a) $Y=Y_o$.}
\smallskip

\emph{b) The $Y$-algebraic space with structural morphism $d_1: Y_1\to Y_o$ is an object of $F(Y)$.}
\smallskip

\emph{c) Over $t$, $Y_{.t}=\mathrm{cosq}_o(f^{-1}(t)/t)$.}

\begin{proof} Notice that for each geometric codimension $1$ point $\overline{y}$ of $Y$ the localization of $S$ at the image $s$ of $\overline{y}$ is noetherian regular of dimension $\leq 1$, since $f$ is by hypothesis flat at $\overline{y}$. If $S$ is local of dimension $1$ with closed point $s$, then with the notations of (10.2) the asserted $S$-groupoid is $\mathrm{cosq}_o(Y/E)$ which, as $Y$ is regular (10.2), is unique up to unique isomorphisms (10.3). By a ``passage \`{a} la limite'' and by gluing, one obtains in the general case an $S$-groupoid $U_.$ satisfying the following conditions :
\smallskip

\emph{a) $U_o$ is open in $Y$ such that $\mathrm{codim}(Y-U_o, Y)\geq 2$.}
\smallskip

\emph{b) The $U_o$-algebraic space with structural morphism $d_1: U_1\to U_o$ is an object of $F(U_o)$.}
\smallskip

\emph{c) Over $t$, $U_{.t}=\mathrm{cosq}_o(f^{-1}(t)/t)$.}
\smallskip

As $d_1: U_1\to U_o$ has the canonical section $s_o: U_o\to U_1$, one may by the proof of (8.2) write $d_1$ as a quotient by a basic finite \'{e}tale $U_o$-equivalence relation on a $U_o$-abelian algebraic space with $s_o$ being the image of the zero section. Hence, by (10.4) there is a cartesian diagram of $S$-algebraic spaces
{\[\xymatrix{
                  U_1 \ar[r]^{} \ar[d]_{d_1}    & Y_1 \ar[d]^{d_1} \\
                  U_o \ar[r]_{}                          & Y_o}\] }where $Y_o=Y$ whose vertical arrow on the right is an object of $F(Y_o)$. By (10.3) this diagram is unique up to unique isomorphisms and there is a unique extension of $s_o: U_o\to U_1$ to a section $s_o: Y_o\to Y_1$ of $d_1: Y_1\to Y_o$. By again (10.3) and arguing as in \cite{almost} 5.7 one finds a unique $S$-groupoid $Y_.$ with $d_1: Y_1\to Y_o$ being a face.

\end{proof}

\smallskip

Similarly as \cite{almost} 3.1, we call every groupoid $Y_.$ satisfying the conditions (10.5) \emph{a})+\emph{b}) an \emph{almost non-degenerate fibration structure} on $f: Y\to S$ with $[Y_.]$ being called the \emph{ramification} $S$-\emph{stack}. Each such fibration has again a tautological factorization 
\[f: Y\to [Y_.]\to S.\]
We say that this structure is \emph{non-degenerate} if the groupoid $Y_.$ is simply connected with $\mathrm{Coker}(d_o, d_1)=S$, that is, if $f: Y\to S$ is an object of $F(S)$ of (10.5). Note that $Y\to [Y_.]$ is then \emph{non-degenerate} in this sense as in \emph{loc.cit.}                   

\smallskip 

{\bf Proposition 10.6.} --- \emph{Keep the notations of $(10.5)$. Consider the following conditions }:
\smallskip

1) \emph{$f$ is proper, $S$ is excellent regular.}
\smallskip

2) \emph{$f$ is proper, $S$ is locally noetherian normal and at each of its points satisfies the condition $(W)$} (EGA IV 21.12.8).
\smallskip

\emph{Then, if $1)$ (resp. $2)$) holds, $S$ is the cokernel of $(d_o, d_1)$ in the full sub-category of the category of $S$-algebraic spaces consisting of the $S$-algebraic spaces (resp. $S$-schemes) which are $S$-separated and locally of finite type over $S$.}

\begin{proof} One argues as in \cite{almost} 5.11.

\end{proof}


\bibliographystyle{amsplain}

\begin{thebibliography}{1}


\bibitem{geometric invariant theory}
\newblock D.~Mumford.
\newblock Geometric Invariant Theory. \emph{Ergebnisse der Mathematik und ihrer Grenzgebiete}, 1965.

\bibitem{raynaud_specialization}
M.~Raynaud.
\newblock Sp\'{e}cialisation du foncteur de Picard.
\newblock \emph{Publications Math\'{e}matiques de l'IH\'{E}S}, 38, 1970.

\bibitem{raynaud_thesis}
M.~Raynaud.
\newblock Faisceaux amples sur les sch\'{e}mas en groupes et les espaces homog\`{e}nes.
\newblock \emph{Lecture Notes in Mathematics}, 119, 1971.

\bibitem{almost}
Y.~Zong.
\newblock Almost non-degenerate abelian fibrations. arxiv.org/abs/1406.5956. 




\end{thebibliography}

\end{document}